\documentclass{article}
\usepackage{amsfonts,amssymb,latexsym,amsmath,amscd,euscript}
\usepackage[all]{xy}
\newcounter{num}[section]
\newcommand{\Num}{ \refstepcounter{num}%
\textbf{\arabic{section}.\arabic{num}}}

\newcommand{\Theorem}{\textbf{Theorem~}}
\newcommand{\Proof}{{\LARGE\emph{Proof}}}
\newcommand{\Def}{\textbf{Definition}}

\newcommand{\Lemma}{ \textbf{Lemma~}}
\newcommand{\Ex}{  \textbf{Example}}
\newcommand{\Remark}{\textbf{Remark}}
\newcommand{\Prop}{\textbf{Proposition~}}
\newcommand{\Cor}{ \textbf{ Corollary~}}

\newcommand{\GL}{{\mathrm{GL}}}

\newcommand{\UT}{{\mathrm{UT}}}

\newcommand{\Db}{{\Bbb D}}
\newcommand{\Pb}{{\Bbb P}}
\newcommand{\Rb}{{\Bbb R}}

\newcommand{\al}{\alpha}
\newcommand{\sal}{{s(\al)}}
\newcommand{\tal}{{t(\al)}}

\renewcommand{\leq}{\leqslant}

\newcommand{\Hc}{{\mathcal{H}}}
\newcommand{\Fc}{{\mathcal{F}}}
\newcommand{\Sc}{{\mathcal{S}}}

\newcommand{\Mat}{{\mathrm{Mat}}}

\newcommand{\Imm}{{\mathrm{Im}}}

\newcommand{\La}{\Lambda}

\begin{document}
\Large

\title{Equidimensional quiver representations and their $U$-invariants}
\author{A.N. Panov}
\date{}
 \maketitle
 {\small   Mechanical and Mathematical  Department, Samara National Research University, Samara, Russia}
  {\small \begin{center}
 		apanov@list.ru
 \end{center}}

 \begin{abstract}
 	
 For an arbitrary equidimensional  quiver representation, we  proposed the method of construction of a system of free generators of the field of $U$-invariants.  The construction of the section and system of generators depends on the choice of a map that assign to each vertex one of the arrows incident to it.

 \end{abstract}
\section{Introduction}
A quiver is a directed graph.  A quiver is given by a pair $Q=(V,A) $, where $V$ is the set of vertices and $A$ is the set of arrows. Each arrow  $\al$ starts at its source 
 $s(\al)$ and terminates at its target $t(\al)$. Multiple arrows and loops are allowed.

Let  $K$ be a field.  A representation of the quiver  $Q$ assigns to each vertex  $v\in V$ a linear space $\{W_v,~ v\in V\}$ defined  over the field  $K$, and to each arrow $\al$ a linear map  $W_{s(\al)}\to W_{t(\al)} $.  Each representation defines the dimensional vector  $n_Q=(n_v)$, where $n_v=\dim W_v$.

We identify each linear space $W_v$  with the coordinate space $K^{n_v}$ and the  space of linear maps  $\Hc_\al$ from  $W_\sal$ to $W_\tal $  with the space of matrices  $\Mat(n_\tal,n_\sal, K)$ of size  
$n_\tal\times n_\sal$.  

 For each vertex $v\in V$, we consider the group  $\GL_v=\GL(n_v)$ with entries in the field  $K$, and its unitriangular group $U_v=\UT(n_v)$ that consists of all upper triangular $(n_v\times n_v)$-matrices  with ones on the diagonal.
Consider the direct product  $\GL_Q=\prod_{v\in V}\GL_v$, its subgroup    $U_Q=\prod_{v\in V} U_v$, and the linear space 
$$\Hc= \Hc_Q=\oplus_{\al\in A} \Hc_\al.$$

The group  $\GL_Q$ acts on $\Hc$ by the formula 
$$g.h=(g_\tal X_\al g_\sal^{-1})_{\al\in A}, \quad g= (g_v) \in \GL_Q, \quad h=(X_\al)\in\Hc.$$
The action of $\GL_Q$ on $\Hc$ determines the representation  $\rho_g$ of the group  $\GL_Q$ in the space of regular functions 
$K[\Hc]$ by the formula $$\rho_gf(h)=f(g^{-1}.h).$$
This representation is extended to the action of $\GL_Q$ on the field of rational functions  $K(\Hc)$.
For any subgroup $G$ of $\GL_Q$ one can define the algebra of invariants $K[\Hc]^G$ and the field of invariants   $K(\Hc)^G$.

 The partial case is the action of the group   $\GL(n)$ on the systems of matrices $g.(X_1,\ldots,X_n)=(gX_1g^{-1},\ldots,gX_ng^{-1})$. This representation is associated  with the quiver  with the set of loops $\al_1,\ldots,\al_n$  with the common vertex. The problem of description of the algebra of $\GL(n)$-invariants is solved in the framework of classical invariant theory in tensors  (see \cite{PV, Pro-1, Don-1}). The algebra of   $\GL(n)$-invariants is generated  by the system of polynomials  $\sigma_t(M)$, where $M$ is an arbitrary monomial of  $X_1,\ldots,X_n$, and  $\{\sigma_t:~~1\leq t\leq n\}$ are coefficients of the characteristic polynomial.  
 
 For an arbitrary quiver $Q$ and a representation of an arbitrary dimension  $(n_v)_{v\in V}$, the algebra of $\GL_Q$-invariants is generated by  the polynomials  $\sigma_t(X_{\beta_1}\cdots X_{\beta_s})$, where  $\beta_1,\ldots,\beta_s$ is an arbitrary   directed closed path  of the quiver.  This theorem in the case of zero characteristic was solved in 1990  in the paper \cite{LBr-Pro} of  L. Le Bruyn, C. Procesi,  and in the case of arbitrary characteristic in 1994 in \cite{Don-2} of S. Donkin.
  
In the present paper, we consider the action of the subgroup $U_Q$ on $\Hc=\Hc_Q$. The problem of construction of a system of generators with their relations
is considered extremely difficult even for the case of coadjoint action of $\UT(n)$ on the space of $\Mat(n)$ (this is associated with the quiver 
of one loop and one vertex).  
 We consider the action of  $U_Q$ on the field $K(\Hc)$. We refer to the field of invariants of this action as the field of  $U$-invariants and  denote  $K(\Hc)^U$. According to the K. Miyata  theorem \cite{Mi}, the field of invariants with respect to an action of  a group of inipotent transformations is rational. Therefore, the field  $K(\Hc)^U$ is rational, i.e. it is a pure transcendental extension of the main field $K$. 
 \\ 
 \Def. We say that a quiver representation is equidimensional if  $n_v=n$ for each $v\in V$. 
  
 In the case of equidimensional representation,  the linear space $\Hc$ is decomposed into a direct sum  of the subspaces $\Hc_\al$, $\al\in A$, which is identified with the space  of   $(n\times n)$-matrices  $\Mat(n)$ defined over the field $K$.
 
  We aim to construct the system of free generators of the field of  $U$-invariants of the equidimensional representation $\Hc$, $n_v=n$,  for an arbitrary quiver  $Q$.  The main result (see Theorem \ref{main-theorem}) is formulated in terms of a fixed map  $\psi:V\to A$ that assign to each vertex $v\in V$ one of the arrows $\psi(v)\in A$  incident to it. The constructed section  $\Sc$ and the system of free generators  $\Pb_\al$,~ $\al\in A$, depend on the choice of the map $\psi$. Earlier in the papers \cite{PanV, Pan-1, Pan-2}, the author constructed such systems for the adjoint representation and representations on the matrix tuples.

We introduce some notations. For positive integers  $a$ and  $b$, we denote by $[a,b]$ the set of all integers $a\leq i\leq b$. 
 For any integer  $1\leq i\leq n$, let  $i'$ stand for the symmetric  number to $i$ with respect to the center of the segment $[1,n]$. We have $i'=n-i+1$.
 
Consider the  $(n\times n)$-matrix $X=(x_{ij})_{i,j=1}^n$ with variables $\{x_{ij}\}$. For each  $1\leq k\leq n$, we denote by  $D_{k}(X)$ the left lower corner minor  of order $k'$ for $X$. In particular, $D_1(X)=\det X$.

 We denote by $\Fc_\al$ the subfield of $K(\Hc_\al)$ generated by  $D_k(X_\al)$,~ $ 1\leq k\leq n$.
 The polynomial  $D_k(X)$ is invariant with respect to the left and right multiplication by  $\UT(n)$.  Therefore for each  $\al\in A$, the polynomials  $\{D_k(X_\al)\}$ are  $U$-invariants, i.e.
  $\{D_k(X_\al)\}$ belongs to $K[\Hc]^U$. 
 Consider the subfield  $\Fc$ of $K(\Hc_\al)^U$ generated by  $\Fc_\al$,~ $\al\in A$. 
 
 We say that  a matrix is anti-diagonal if its entries above and below the anti-diagonal are zeros. We denote by  $\La(n)$ the subspace of all anti-diagonal  $(n\times n)$-matrices over the field  $K$. We consider the subspace of upper anti-triangular matrices  $\Sc^+(n)$ over the field  $K$
(their entries below the anti-diagonal are zeros) and the subspace of lower anti-triangular matrices $\Sc^-(n)$ (their entries above the anti-diagonal are zeros). 

We introduce a linear order on the set of all pairs 
$\{(a,b)\}$, ~ $1\leq a,b \leq n$, as follows: $(a,b)\prec (i,k)$, if $b<k$ or $b=k$ and $a>i$.  Then $(n,1)\prec (n-1,1)\prec \ldots\prec (1,1)\prec (n,2)\prec\ldots\prec (1,n)$.
\\
\Def. Let $X=(x_{ij})_{i,j=1}^n$ be the matrix with variables $\{x_{ij}\}$. \\
~1) We say that a polynomial  $f(X)$ is obtained by a triangular transformation from the variable $x_{ik}$ if  $f(X)=f_0(X)x_{ik}+f_1(X)$, where $f_0(X)$ and $f_1(X)$ are polynomials in variables  $x_{ab}$,~ $(a,b)\prec (i,k)$ defined over the same field  as $f(X)$, and  $f_0(X)\ne 0$.\\
2) We say that the system of polynomials  $\{f_{ik}(X):~ 1\leq i, k\leq n\}$, is obtained by a triangular transformation from the system  $\{x_{ik}:~ 1\leq i, k\leq n\}$, if each $f_{ik}(X)$ is obtained by a triangular transformation from  $x_{ik}$. \\
\Remark. If $\{f_{ik}(X):~ 1\leq i, k\leq n\}$ is obtained by a a triangular transformation from  $\{x_{ik}:~ 1\leq i, k\leq n\}$, then the first system of polynomials (and the second also) freely generate the field of rational functions  $K(X)$.

\section{Construction of section}

 In this section, we construct a section for the action of the subgroup  $U_Q$ on $\Hc$.   
 Fix a map  $\psi:V\to A$ such that for each  $v\in V$ the arrow $\psi(v)$ is incident to  $v$. Then for any $\al$ from $\Imm\,\psi$ we have $\al=\psi(\sal)$ or $\al=\psi(\tal)$. It is possible that  $\al=\psi(\tal)=\psi(\sal)$.    
 The construction of the section and the system of free generators depends on the choice of  the map $\psi$.
 
We consider the  subspace  $\Sc=\oplus_{\al\in A} \Sc_\al$, where the subspace $\Sc_\al$ in $\Hc_\al$ is defined as follows. If $\al$  connects different vertices  (i.e. $\xymatrix{\tal&\sal\ar[l]_\al}$ and $\sal\ne\tal$),  we define 

\begin{equation}\label{Sec-diff}
	\Sc_\al = \left\{ \begin{array}{l}
		\Hc_\al= \Mat(n),\mbox{if}~\al\notin \Imm\,\psi,\\
		\Sc^-(n),\mbox{if}~\al =\psi(\tal)\ne\psi(\sal),\\
		\Sc^+(n),\mbox{if}~\al =\psi(\sal)\ne\psi(\tal),\\
		\La(n), \mbox{if}~\al =\psi(\tal)=\psi(\sal).
	\end{array}
	\right.
\end{equation}

If $\al$ is a loop, we denote $q=\sal=\tal$. We have $\xymatrix{q\ar@(ul,dl)[]_\al}$. We define
\begin{equation}\label{Sec-loop}
	\Sc_\al = \left\{ \begin{array}{l}
		\Hc_\al= \Mat(n),~\mbox{if}~\al\ne\psi(q),\\
		\Sc^-(n),~ \mbox{if}~\al =\psi(q).\\
	\end{array}
	\right.
\end{equation}
\Prop\Num\label{section}. The space  $\Hc$ is a closure of the set   $$\bigcup_{g\in U_Q}g.\Sc .$$
\Proof.
Consider the open subset  $\Omega\subset \Hc$ that consists of all  $h=(X_\al)_{\al\in A}$ such that  $D_k(X_\al)\ne 0$ for all $1\leq k\leq n$. Let us show that for any  $h=(X_\al)_{\al\in A}\in\Omega$ there exist  $g=(g_v)_{v\in V}\in U$ such that  $g.h\in \Sc$, i.e. 
\begin{equation}\label{in-section-al}
(g.h)_\al = g_\tal X_\al g_\sal^{-1}\in \Sc_\al \quad\mbox{for~~each}\quad \al\in A.	
\end{equation}

If $\al\notin \Imm\,\psi$, then $\Sc_\al=\Mat(n)$, and  (\ref{in-section-al}) is true for all  $g=g=(g_v)_{v\in V}$.

One knows the following property: for any $(n\times n)$-matrix  $X$ obeying the condition  $D_k(X)\ne 0$,~ $1\leq k\leq n$, there exist 
matrices $u_1, u_2\in \UT(n)$ such that   $u_1 X u^{-1}_2\in\La(n)$ (it implies  $u_1X\in \Sc^-(n)$  and  $Xu^{-1}_2\in \Sc^+(n)$). 

Consider the case when   $\al$ connects different vertices. 
If $\al =\psi(\tal)\ne\psi(\sal)$, then we have $\Sc_\al=\Sc^-(n)$. Taking $X=X_\al$ and  $g_\tal=u_1$, we obtain  $g_\tal X_\al\in \Sc^-(n)$. Then  $g_\tal X_\al g_\sal^{-1}$ belongs to  $ \Sc^-(n)$ for any $g_\sal$.  We get $(g.h)_\al\in \Sc_\al$.

If $\al =\psi(\sal)\ne\psi(\tal)$, then we have $\Sc_\al=\Sc^+(n)$. Taking $g_\sal=u_2$, we obtain  $ X_\al g_\sal^{-1}\in \Sc^+(n)$. Then  $g_\tal X_\al g_\sal^{-1}$ belongs to $ \Sc^+(n)$ for any $g_\tal$.  We get $(g.h)_\al\in \Sc_\al$.

If $\al =\psi(\tal)=\psi(\sal)$, then we have $\Sc_\al=\La(n)$. Taking  $g_\tal=u_1$ and $g_\sal=u_2$, we obtain that   $g_\tal X_\al g_\sal^{-1}$ belongs to $ \La(n)$.  We get $(g.h)_\al\in \Sc_\al$.

Consider the case when $\al$ is a loop with vertex   $q$. If $\al=\psi(q)$, then we have $\Sc_\al=\Sc^-(n)$. Taking $g_q=u_1$, we obtain $g_q X_\al\in \Sc^-(n)$. Then  $g_q X_\al g_q^{-1}$ belongs to $ \Sc^-(n)$.  We get $(g.h)_\al\in \Sc_\al$. ~ $\Box$

Proposition \ref{section} impies the  corollary.\\
\Cor\Num\label{restriction}.
The restriction map  $\pi: K[\Hc]^U\to K[\Sc]$ is an embedding.
The embedding  $\pi$ extends to the embedding of fields  $$\pi: K(\Hc)^U \hookrightarrow K(\Sc).~~\Box$$

In particular, the field $\pi(\Fc_\al)$ is a subfield of $K(\Sc_\al)$. If $\Sc_\al$ equals to $\Sc^\pm(n)$ or $\La(n)$, then $\pi(\Fc_\al)$ is generated by the anti-diagonal elements of the matrix $X_\al\in \Sc_\al$.

\section{$U$-invariants associated with arrows}

In this section, we assign to each arrow $\al\in A$ a system of polynomials $\Pb_\al$. These polynomials are $U$-invariants.
We show that the system  $\pi(\Pb_\al\sqcup \Db_\al)$ is obtained by a triangular transformation over the field $\pi(\Fc)$  from the system of matrix entries $\{x_{ab}\} $ of the matrix $X_\al$,  which form  the standard basis in the dual space of $\Sc_\al$. Moreover,  $\pi( \Db_\al)$ is obtained by a triangular transformation over the field  $K$ from the matrix entries  $\{x_{a,a'}\} $ that belong to the anti-diagonal of $X_\al$, and $\pi( \Pb_\al)$ is obtained by a triangular transformation over the field  $\pi(\Fc)$  from the other matrix entries of $X_\al$. 

If $\al$ connects different vertices $\xymatrix{\tal&\sal\ar[l]_\al}$, then we define  $\Pb_\al=\Pb_\sal\sqcup\Pb_\tal$, where $\Pb_\sal$ (respectively,  $\Pb_\tal)$ is a system of   $U$-invariants associated with the source $\sal$ (respectively, with the target $\tal$) of the arrow $\al$ (see formula  (\ref{Pb-alpha})). 

If $\Sc_\al=\Sc^+(n)$ or $\Sc_\al=\La(n)$, then we take $\Pb_\sal=\varnothing$. If  $\Sc_\al=\Mat(n)$ or $\Sc_\al=\Sc^-(n)$, the system of standard coordinate functions  $\{x_{ab}: a'<b\} $ (they lie below the anti-diagonal of $X_\al$) is algebraically independent on $\Sc_\al$.  
The system $\Pb_\sal$  has the property that its restriction  $\pi(\Pb_\sal)$ on the section $\Sc$ is obtained by a triangular transformation over $\pi(\Fc)$ from  $\{x_{ab}: a'<b\} $ (see Lemmas \ref{lem-one} and \ref{lem-two}).

Similarly, we define $\Pb_\tal$.  If $\Sc_\al=\Sc^-(n)$ or $\Sc_\al=\La(n)$, then we take $\Pb_\tal=\varnothing$.  If $\Sc_\al=\Mat(n)$ or $\Sc_\al=\Sc^+n)$, then the system of standard coordinate functions  $\{x_{ab}: a'>b\} $ (they lie above the anti-diagonal of $X_\al$) is algebraically independent on $\Sc_\al$. 
The system $\Pb_\tal$  has the property that its restriction  $\pi(\Pb_\tal)$ on the section $\Sc$ is obtained by a triangular transformation over $\pi(\Fc)$ from   $\{x_{ab}: a'>b\} $ (see Lemmas \ref{lem-three} and \ref{lem-four}). 

If $\al$ is a loop, then we take  $\Pb_\tal=\varnothing$ and $\Pb_\al=\Pb_\sal$.

In what follows, in Cases 1 and 2, for each arrow $\beta$ incident to its source $\sal$, we construct a system of  $U$-invariants  $\Pb_{\sal,\beta}$. We define  $\Pb_\sal=\Pb_{\sal,\beta}$, where $\beta=\psi(\sal)$   (see  (\ref{Pb-sal})). 

Respectively, in Cases  3 and 4, for each  $\gamma$ incident to its target $\tal$, we construct a system of  $U$-invariants  $\Pb_{\gamma,\tal}$. We define $\Pb_\tal=\Pb_{\gamma,\tal}$, where $\gamma=\psi(\tal)$    (see (\ref{Pb-tal})). 
\\
\\ 
{\bf Case 1}. Let  $\xymatrix{\tal&\sal\ar[l]_\al &v\ar[l]_\beta}$, i.e. $\beta$ terminates at $\sal$. We admit that  $\al$  or  $\beta$ is a loop (i.e. $\tal=\sal$ or $\sal=v$).  It is possible $\al$ and $\beta$ are loops and $\al=\beta$.

The action of $g\in \GL_Q$ on the pair of matrices $(X,Y)$, where $X=X_\al$ and $Y=X_\beta$, is defined by the formula
$g.(X,Y) = (g_\tal Xg_\sal^{-1}, g_\sal Yg_v^{-1})$.  

For the pair $i'\leq j$ (i.e.  $(i,j)$ lies on or below the anti-diagonal), we denote by  $M_{ij}(X)$ the minor of order  $i'$ of the matrix  $X$ with 
the system of rows  $[i,n]$ and columns  $[1,i'-1]\sqcup \{j\}$.

For the pair $j\leq k$, we denote by $N_{jk}(Y)$ the minor of order  $k'$ of the matrix   $Y$ with the system of rows $\{j\}\sqcup [k+1,n]$ and columns $[1,k']$. 
Observe that  $$D_k(X) = M_{k,k'}(X),\qquad D_k(Y) = N_{k,k}(Y).$$  

Let $i'<k$ (i.e.  $(i,k)$ lies below the anty-diagonal). We define  the polynomial  
\begin{equation}\label{PXY}
P_{ik}(X,Y) = \sum_{i'\leq j\leq k} M_{ij}(X) N_{jk}(Y).
\end{equation}

The polynomials $$\Pb(X,Y)=\{P_{ik}(X,Y):~i'<k\}$$ are  $U$-invariants \cite{Pan-2}.\\ 

\Lemma\Num\label{lem-one}. 1)  The polynomial $P_{ik}(X,Y)$, where $(i,k)$ lies below the anti-diagonal, is obtained by a triangular transformation from  $x_{ik}$ over the field  $K(Y)$.
\\
2)  Let $\beta$ obey the condition of  Case 1, and  $\beta=\psi(\sal)$. Then the system of standard coordinate functions  $\{x_{ab}: a'<b\} $ (they lie below the anti-diagonal of $X_\al$) is algebraically independent on $\Sc_\al$.  
The system $\pi(\Pb(X_\al, X_\beta)) $  is obtained by a triangular transformation over $\pi(\Fc)$ from   $\{x_{ab}: a'<b\} $. \\
\\
\Proof.  Statement 1) follows from the fact that each minor $M_{ij}(X)$,~$i'\leq j\leq k$, from formula (\ref{PXY}) is obtained by a triangular transformation from  $x_{ij}$. 

Let   $\beta$ obey the conditions of item  2). Then    either  $\Sc_\al$ equals to   
$\Mat(n)$ (if $\al\ne\psi(\tal)$, or $\Sc^-(n)$ (if $\al=\psi(\tal)$). 
The system of standard coordinate functions  $\{x_{ab}: a'<b\} $  is algebraically independent on $\Sc_\al$.  .

From  $\beta=\psi(\sal)$, it implies  that  either  $\Sc_\beta=\Sc^-(n)$, or  $\Sc_\beta=\La(n)$.
Let $Y=S^-$  be the matrix with entries which are the standard coordinate functions on $ \Sc_\beta$. 
Then  $N_{jk}(S^-)=0$ for each $i'\leq j< k$. Therefore, $$P_{ik}(X,S^-)= M_{ik}(X) N_{k,k}(S^-) = M_{ik}(X) D_k(S^-).$$
 The polynomial $P_{ik}(X,S^-)$ is obtained by a triangular transformation over $\pi(\Fc_\beta)$ from  $x_{ik}$ , and the system of polynomials $\Pb(X,S^-)$
is obtained by a triangular transformation over the field $\pi(\Fc_\beta)$ from  $\{x_{ab}\}$ lying below the anti-diagonal of $X$. This proves 2). ~$\Box$\\
\\
{\bf Example 1} 
{\small $$ X= \left(\begin{array}{ccc}
x_{11}&x_{12}&x_{13}\\
x_{21}&x_{22}&x_{23}\\
x_{31}&x_{32}&x_{33}
\end{array}\right), \qquad Y= \left(\begin{array}{ccc}
y_{11}&y_{12}&y_{13}\\
y_{21}&y_{22}&y_{23}\\
y_{31}&y_{32}&y_{33}
\end{array}\right),\qquad S^-= \left(\begin{array}{ccc}
0&0&s_{13}\\
0&s_{22}&s_{23}\\
s_{31}&s_{32}&s_{33}
\end{array}\right),$$
 $$D_3(X) = x_{31},\quad D_2(X)= \left|\begin{array}{cc}
x_{21}&x_{22}\\
x_{31}&x_{32}
\end{array}\right|, \quad D_1(X)=\det(X),$$
$$D_3(Y) = y_{31},\quad D_2(Y)= \left|\begin{array}{cc}
y_{21}&y_{22}\\
y_{31}&y_{32}
\end{array}\right|, \quad D_1(Y)=\det(Y),$$}
{ \small 
	$$P_{32}(X,Y) = 
	x_{31} \left|\begin{array}{cc}
	y_{11}&y_{12}\\
	y_{31}&y_{32}
	\end{array}\right| + x_{32} \left|\begin{array}{cc}
	y_{21}&y_{22}\\
	y_{31}&y_{32}
	\end{array}\right|,\qquad P_{32}(X,S^-) = 
	x_{32} \left|\begin{array}{cc}
	0&s_{22}\\
	s_{31}&s_{32}
	\end{array}\right| = D_2(S^-) x_{32},$$
	
	$$P_{33}(X,Y) = 
	x_{31}y_{11}+x_{32}y_{21}+x_{33}y_{31},\qquad P_{33}(X,S^-) = 
	x_{33}s_{31}=D_3(S^-)x_{33},$$
	
	$$P_{23}(X,Y) = 
	\left|\begin{array}{cc}
	x_{21}&x_{22}\\
	x_{31}&x_{32}
	\end{array}\right| y_{21} + \left|\begin{array}{cc}
	x_{21}&x_{23}\\
	x_{31}&x_{33}
	\end{array}\right| y_{31}, \qquad P_{23}(X,S^-) = 
	\left|\begin{array}{cc}
	x_{21}&x_{23}\\
	x_{31}&x_{33}
	\end{array}\right| s_{31} = D_3(S^-) 	\left|\begin{array}{cc}
	x_{21}&x_{23}\\
	x_{31}&x_{33}
	\end{array}\right|.$$}
{\bf Case 2}. Let
  $\xymatrix{\tal&\sal\ar[l]_\al \ar[r]^\beta & v }$, i.e.  $\beta$ starts at  $\sal$. Here we assume that $\al\ne \beta$, the arrow $\al$ may be a loop, and  $\beta $ is not a loop (the case $\beta$ is a loop,  is considered in Case  1).

The action of $g\in \GL_Q$ on the pair of matrices $(X,Y)$, where $X=X_\al $ and $Y=X_\beta$,  is defined by the formula 
$g.(X,Y) = (g_\tal Xg_\sal^{-1}, g_v Y g_\sal^{-1})$.

Consider the $(2n\times n)$-matrix $\left(\frac{X}{Y}\right)$. Let the pair  $(i,k)$ lie below the anti-diagonal (i.e. $i'<k$). Decompose $k$ into a sum of positive integers  $k= i'+(k-i')$.   
Consider the  $(k\times k)$-matrix $$\left(\frac{X_{ik}}{Y_{ik}}\right),$$
where  $X_{ik}$  is the left lower block of the matrix  $X$  of size  $i'\times k$, and $Y_{ik}$ is the left  lower  block of the matrix  $Y$ of size $(k-i')\times k$. 
Let us show that the determinant of this matrix 
$$R_{ik}^-(X,Y) =\det \left(\frac{X_{ik}}{Y_{ik}}\right)$$
is a  $U$-invariant. 
For an arbitrary matrix  $A$, let  $A_t$ stand for its right lower   $(t\times  t)$-block,
and $A^{(t)}$ is its left upper   $(t\times t)$-block.

In these notations 
$$g.\left(\frac{X_{ik}}{Y_{ik}}\right) = \left(\begin{array}{cc}
\left(g_\tal\right)_{i'}&0\\
0&\left(g_v\right)_{k-i'}
\end{array}\right) \left(\frac{X_{ik}}{Y_{ik}}\right) \left(g^{-1}_\sal\right)^{(k)}.
$$
Since the matrices  $g_\tal$, $g_\sal$, $g_v$  belong to $\UT(n)$, we have $$R_{ik}^-(g_\tal Xg_\sal^{-1}, g_v Y g_\sal^{-1}) = R_{ik}^-(X,Y).$$ The polynomials  $\{R_{ik}^-(X,Y)\}$ are $U$-invariants. 

Introduce the system of $U$-invariants
$$\Rb^-(X,Y) =\{R_{ik}^-(X,Y): ~i'<k\}.$$
\Lemma\Num\label{lem-two}.  1) The polynomial $R^-_{ik}(X,Y)$, where $(i,k)$ lies below the anti-diagonal, is obtained by a triangular transformation over the field $K(Y)$ from  $x_{ik}$. \\
2)  Let $\beta$ obey the condition of Case 2,   and  $\beta=\psi(\sal)$. 
Then the system of standard coordinate functions  $\{x_{ab}: a'<b\} $ (they lie below the anti-diagonal of $X_\al$) is algebraically independent on $\Sc_\al$.  
The system $\pi(\Rb^-(X_\al,X_\beta)) $  is obtained by a triangular transformation over $\pi(\Fc)$ from   $\{x_{ab}: a'<b\} $. \\
\\
\Proof.  The statement 1) follows from the fact that the element  $x_{ik}$ lies in the right upper corner of the determinant   $R^-_{ik}(X,Y)$. 

 Let $\beta$ obey the conditions of item  2).  Since  $\beta\ne\al$ and $\beta=\psi(\sal)$, we have $\al\ne\psi(\sal)$. Therefore,  $\Sc_\al$ is either 
 $\Mat(n)$, or  $\Sc^-(n)$.
 The system of standard coordinate functions  $\{x_{ab}: a'<b\} $ is algebraically independent on  $\Sc_\al$.

It follows from  $\beta=\psi(\sal)$ that  either
    $\Sc_\beta=\Sc^+(n)$, or $\Sc_\beta=\La(n)$.
Let $Y=S^+$ be the matrix with entries, which are the standard coordinate functions on $\Sc_\beta$.  Then  $R^-_{ik}(X,S^+)$ decomposes into a product 
$$R^-_{ik}(X,S^+) = \large{\left| \begin{array}{cccccc}
	x_{i1}&\ldots&x_{i,k-i'}&x_{i,k-i'+1}&\ldots&x_{ik}\\
	\cdots&\cdots&\cdots&\cdots&\cdots&\cdots\\
		x_{n1}&\ldots&x_{n,k-i'}&x_{n,k-i'+1}&\ldots&x_{nk}\\
			y_{i'+k',1}&\ldots&	y_{i'+k',k-i'}&	0&\ldots&0\\
				\cdots&\cdots&\cdots&\cdots&\cdots&\cdots\\
					y_{n1}&\ldots&0&0&\ldots&0\\
\end{array}
\right|}
=
\pm M_I^J(X) D_{k'+i'}(S^+),$$
where $M_I^J(X)$ is the minor of the matrix  $X$ with the system of rows  $I=[i,n]$ and columns  $J=[k-i'+1,k]$.
The element $x_{ik}$ lies in the right upper corner of the minor. Therefore, 
$R^-_{ik}(X,S^+)$  is obtained by a triangular transformation over the field  $\pi(\Fc_\beta)$ from    $x_{ik}$ , and the system of polynomials  $\Rb^-(X,S^+)$ is obtained by a triangular transformation over the field $\pi(\Fc_\beta)$ from the system of coordinate functions  $\{x_{ab}\}$ lying below the anti-diagonal of  $X$.
 This proves  2). ~$\Box$
 \\
 {\bf Example 2}. For  $n=3$ we have  
{ \small $$ \left(\frac{X}{Y}\right) =  \left(\begin{array}{ccc}
	x_{11}&x_{12}&x_{13}\\
	x_{21}&x_{22}&x_{23}\\
	x_{31}&x_{32}&x_{33}\\
		\hline
		y_{11}&y_{12}&y_{13}\\
	y_{21}&y_{22}&y_{23}\\
	y_{31}&y_{32}&y_{33}\\
	\end{array}
	\right),\qquad  \left(\frac{X}{S^+}\right) =  \left(\begin{array}{ccc}
	x_{11}&x_{12}&x_{13}\\
	x_{21}&x_{22}&x_{23}\\
	x_{31}&x_{32}&x_{33}\\
	\hline
	s_{11}&s_{12}&s_{13}\\
	s_{21}&s_{22}&0\\
	s_{31}&0&0
	\end{array}
	\right),$$}

{ \small
	$$R^-_{32} (X,Y) = 
	\left|\begin{array}{cc}
	x_{11}&x_{32}\\
		y_{31}&y_{32}\\
	\end{array}\right|,\qquad R^-_{32} (X,S^+) = 
	\left|\begin{array}{cc}
	x_{31}&x_{32}\\
	s_{31}&0\\
	\end{array}\right| =-D_3(S^+)x_{32},$$
	$$R^-_{33} (X,Y) = 
	\left|\begin{array}{ccc}
	x_{31}&	x_{32}&x_{33}\\
	y_{21}&y_{22}&y_{23}\\
	y_{31}&	y_{32}&y_{33}
	\end{array}\right|,\qquad R^-_{33} (X,S^+) = 
	\left|\begin{array}{ccc}
	x_{31}&	x_{32}&x_{33}\\
	s_{21}&s_{22}&0\\
	s_{31}&	0&0
	\end{array}\right|=  D_2(S^+) x_{33}$$,	
	$$ 	
	R^-_{23} (X,Y) = 
	\left|\begin{array}{ccc}
	x_{21}&	x_{22}&x_{23}\\
	x_{31}&	x_{32}&x_{33}\\
	y_{31}&	y_{32}&y_{33}
	\end{array}\right|,\qquad R^-_{23} (X,S^+) = 
	\left|\begin{array}{ccc}
		x_{21}&	x_{22}&x_{23}\\
	x_{31}&	x_{32}&x_{33}\\
	s_{31}&	0&0
	\end{array}\right|=  D_3(S^+) \left|\begin{array}{cc}
	x_{22}&x_{23}\\
	x_{32}&x_{33}\\
	\end{array}\right|.		 $$ }

We define 
 $$
\Pb_{\sal,\beta} (X_\al,X_\beta)=\left\{\begin{array}{l} \Pb(X_\al, X_\beta), \mbox{if}~\beta~\mbox{obeys~ the~ conditions~ of~ Case~1},\\	
	 \Rb^-(X_\al, X_\beta), \mbox{if}~\beta~\mbox{obeys~ the~ conditions~ of~ Case~2}.
\end{array}\right.
$$
For $h=(X_\al)\in\Hc$, we take
\begin{equation}\label{Pb-sal}
	\Pb_{\sal}(h) = \left\{\begin{array}{l} \varnothing, \quad\mbox{if}\quad \al=\psi(\sal),\\
			\Pb_{\sal, \beta} (X_\al, X_\beta),~ \mbox{where}~ \beta = \psi(\sal), \quad\mbox{if}\quad \al\ne\psi(\sal). \end{array}\right.
\end{equation}
Let us formulate a statement that follows from Lemmas  \ref{lem-one} and   \ref{lem-two}\\
\Prop\Num\label{prop-sal}. Let $\al$ obeys the conditions of Case 1 or 2, and $\beta=\psi(\sal)$. Then the system of standard coordinate functions  $\{x_{ab}: a'<b\} $ (they lie below the anti-diagonal of $X_\al$) is algebraically independent on $\Sc_\al$.  
The system  of polynomials $\pi(\Pb_\sal(X_\al,X_\beta)) $ is obtained by a triangular transformation over $\pi(\Fc)$ from   $\{x_{ab}: a'<b\} $.

Let us turn to a construction of the system of polynomials  $\Pb_\tal$ associated with the target  $\tal$ of the arrow $\al$. Observe that if $\al$ is a loop (i.e. $\sal=\tal=q$), then the pair $\al$ and $\beta=\psi(q)$ satisfies conditions of Case 1 or Case 2. In Cases 3 and 4, $\al$ is not a loop.\\
\\
 {\bf Case 3}. Let $\xymatrix{v&\tal\ar[l]_\gamma &\sal\ar[l]_\al}$, i.e.  $\gamma$ starts at  $\tal$. Here we assume that  $\al$ and $\gamma$ are not loops (i.e. $\tal\ne \sal$ and $v\ne \tal$).

The action of  $g\in \GL_Q$ on the pair of matrices $(Z,X)$, where $Z=X_\gamma$ and $X=X_\al$, is defined by the formula 
 $g.(Z,X) = (g_vZg_\tal^{-1}, g_\tal Xg_\sal^{-1})$.
The polynomials $$\Pb(Z,X)=\{P_{ik}(Z,X):~i'<k\}$$ are  $U$-invariants (see Case 1 and \cite{Pan-2}).\\
\\
\Lemma\Num\label{lem-three}.
1) The polynomial $P_{ik}(Z,X)$, where $(i,k)$ lies below the anti-diagonal, is obtained by a triangular transformation from $x_{i',k'}$ over  $K(Z)$ (observe that if  $(i,k)$ lies below the anti-diagonal, then $(i',k')$ lies above the anti-diagonal). \\
2) Let $\gamma$ obey the conditions of Case 3, and  $\gamma=\psi(\tal)$. Then the system of standard coordinate functions  $\{x_{ab}: a'>b\} $ (they lie above the anti-diagonal of $X_\al$) is algebraically independent on $\Sc_\al$.  
The system $\pi(\Pb(X_\gamma, X_\al)) $  is obtained by a triangular transformation over $\pi(\Fc_\gamma)$ from   $\{x_{ab}: a'>b\} $. \\
\\
\Proof. Statement 1)  follows from the fact that each minor  $N_{jk}(X)$,~$i'\leq j\leq k$  from formula (\ref{PXY}), where $i'<k$, is obtained by a triangular transformation from  $x_{j,k'}$, $ j\leq k$. 

Let  $\gamma$ obey the conditions of item  2). The condition  $\gamma=\psi(\tal)$ implies $\al\ne\psi(\tal)$. Therefore, $\Sc_\al$ equals to 
$\Mat(n)$ (if $\al\ne\psi(\sal)$, or $\Sc^+(n)$ (if $\al=\psi(\sal)$).
The system of standard coordinate functions  $\{x_{ab}: a'>b\} $  is algebraically independent on $\Sc_\al$.  
 
It follows from  $\gamma=\psi(\tal)$ that  either  $\Sc_\gamma=\Sc^+(n)$ or 
$\Sc_\gamma=\La(n)$.
Let  $Z=S^+$  be the matrix with entries, which are the standard coordinate functions on $\Sc_\gamma$. Then  $M_{jk}(S^+)=0$ for each $i'< j\leq k$. Therefore,  $$P_{ik}(S^+,X)= N_{i,i'}(S^+) N_{i',k}(X) = D_i(S^+) N_{i',k}(X).$$
 The polynomial  $P_{ik}(S^+,X)$  is obtained by a triangular transformation over $\pi(\Fc_\gamma)$ from  $x_{i',k'}$, and the system of polynomials  $\Pb(S^+,X)$ is obtained by a triangular transformation over $\pi(\Fc_\gamma)$ from  $\{x_{ab}\}$ lying above the anti-diagonal of $X$.
This proves 2). ~$\Box$
\\
{\bf Example 3}. For $n=3$ we have  
{ \small $$ Z= \left(\begin{array}{ccc}
	z_{11}&z_{12}&z_{13}\\
	z_{21}&z_{22}&z_{23}\\
	z_{31}&z_{32}&z_{33}
	\end{array}\right)\qquad X= \left(\begin{array}{ccc}
	x_{11}&x_{12}&x_{13}\\
	x_{21}&x_{22}&x_{23}\\
	x_{31}&x_{32}&x_{33}
	\end{array}\right), \qquad ,\qquad S^+= \left(\begin{array}{ccc}
	s_{11}&s_{12}&s_{13}\\
	s_{21}&s_{22}&0\\
	s_{31}&0&0
	\end{array}\right),$$}
{ \small 
	$$P_{32}(Z,X) = 
	z_{31} \left|\begin{array}{cc}
	x_{11}&x_{12}\\
	x_{31}&x_{32}
	\end{array}\right| + z_{32} \left|\begin{array}{cc}
	x_{21}&x_{22}\\
	x_{31}&x_{32}
	\end{array}\right|,\qquad 
	P_{32}(S^+,X) = 
	s_{31} \left|\begin{array}{cc}
	x_{11}&x_{12}\\
	x_{31}&x_{32}
	\end{array}\right| = D_3(S^+) \left|\begin{array}{cc}
	x_{11}&x_{12}\\
	x_{31}&x_{32}
	\end{array}\right|,$$
		
	$$P_{33}(Z,X) = 
	z_{31}x_{11}+z_{32}x_{21}+z_{33}x_{31},\qquad P_{33}(S^+,X) = 
	s_{31}x_{11}=D_3(S^+)x_{11},$$
	
		$$P_{23}(Z,X) = 
	\left|\begin{array}{cc}
	z_{21}&z_{22}\\
	z_{31}&z_{32}
	\end{array}\right| x_{21} + \left|\begin{array}{cc}
	z_{21}&z_{23}\\
	z_{31}&z_{33}
	\end{array}\right| x_{31}, \qquad P_{23}(S^+,X)=	\left|\begin{array}{cc}
	s_{21}&s_{22}\\
	s_{31}&0
	\end{array}\right| x_{21}= D_2(S^+) x_{21}.$$}
{\bf Case 4}. Let  $\xymatrix{v\ar[r]^\gamma&\tal &\sal\ar[l]_\al}$, i.e. $\gamma$ terminates at  $\tal$. Here we assume that $\al\ne\gamma$,  $\al$ is not a loop, 
and $\gamma$ maybe a loop.

The action of  $g\in \GL_Q$ on the pair of matrices  $(Z,X)$, where $Z=X_\gamma$  and $X=X_\al$, is defined by the formula 
$g.(Z,X) = (g_\tal Zg_v^{-1}, g_\tal Xg_\sal^{-1})$.

Consider the  $(n\times 2n)$-matrix  $(Z|X)$. Let the pair  $(i,k)$ lie above the anti-diagonal (i.e. $i'>k$). Decompose $i'$ into a sum of positive integers  $i'= k+(i'-k)$. Consider the $(i'\times i')$-matrix $ (Z_{ik}|X_{ik})$,
where $Z_{ik}$ is the left lower block of the matrix  $Z$ of size $i'\times(i'-k)$, and $X_{ik}$ is the  left lower block of the matrix  $X$ of size $i'\times k$.

Let us show that the determinant of this matrix  
$$R_{ik}^+(Z,X) =\det (Z_{ik}|X_{ik})$$
is  $U$-invariant. 
In notations we used in Case 2, we have 
$$g. (Z_{ik}|X_{ik})= \left(g_\tal\right)_{i'} (Z_{ik}|X_{ik})\left(\begin{array}{cc}
	\left(g_v^{-1}\right)^{(i'-k)}&0\\
	0&\left(g^{-1}_\sal\right)^{(k)}
\end{array}\right).
$$
Since the matrices  $g_\tal$, $g_\sal$, $g_v$ belong to   $\UT(n)$, we have $$R_{ik}^+(g_\tal Zg_v^{-1}, g_\tal Xg_\sal^{-1}) = R_{ik}^+(Z,X).$$ The polynomials $\{R_{ik}^+(X,Y)\}$ are $U$-invariants. 

We introduce the system of $U$-invariants
$$\Rb^+(Z,X) =\{R_{ik}^+(Z,X): ~i'>k\}.$$
\Lemma\Num\label{lem-four}.  1) The polynomial $R^+_{ik}(Z,X)$, where $(i,k)$ lies above the anti-diagonal, is obtained by a triangular transformation from  $x_{ik}$ over the field  $K(Z)$.\\
2)  Let $\gamma$ obey the conditions of Case 4,   and $\gamma=\psi(\tal)$. Then the system of standard coordinate functions  $\{x_{ab}: a'>b\} $ (they lie above the anti-diagonal of $X_\al$) is algebraically independent on $\Sc_\al$.  
The system  $\pi(\Rb^+(X_\gamma, X_\al)) $  is obtained by a triangular transformation over $\pi(\Fc)$ from   $\{x_{ab}: a'>b\} $. 
\\
\\
\Proof. The statement 1) follows from the fact that the element $x_{ik}$ lies in the right upper corner of the matrix   $R^+_{ik}(Z,X)$. 

 Let $\gamma$ obey the conditions of item 2). It follow from  $\gamma\ne\al$ and $\gamma=\psi(\tal)$ that 
$\al\ne\psi(\tal)$.  Therefore $\Sc_\al$ equals to either
$\Mat(n)$ (if $\al\ne\psi(\sal)$, or $\Sc^+(n)$ (if $\al=\psi(\sal)$).  The system of matrix standard coordinate functions $\{x_{ab}\}$ lying above the anti-diagonal is  algebraically on $\Sc_\al$.

From  $\gamma=\psi(\tal)$, it implies that either  $\Sc_\gamma=\Sc^-(n)$, or 
$\Sc_\gamma=\La(n)$.
Let  $Z=S^-$  be the matrix with entries, which are the standard coordinate functions on $\Sc_\gamma$.
Then   $R^-_{ik}(S^-,X)$ decomposes into a product  
$$R^+_{ik}(S^-,X) = \large{\left| \begin{array}{cccccc}
		0&\ldots&0&x_{i1}&\ldots&x_{ik}\\
		\cdots&\cdots&\cdots&\cdots&\cdots&\cdots\\
		0&\ldots&0&x_{i+k-1,1}&\ldots&x_{i+k-1,k}\\
		0&\ldots& z_{i+k,i'-k}&	x_{i+k,1}&\ldots&	x_{i+k,k}\\
		\cdots&\cdots&\cdots&\cdots&\cdots&\cdots\\
		z_{n1}&\ldots&z_{n,i'-k}&x_{n1}&\ldots&x_{nk}\\
	\end{array}
	\right|}
=
 \pm D_{i+k}(S^-)  N_I^J(X),$$
where $N_I^J(X)$is the minor of the matrix  $X$ with the system of rows  $I=[i,i+k-1]$ and columns  $J=[1,k]$.
The element $x_{ik}$ lies in the right upper cornor of this minor.  Therefore, 
$R^+_{ik}(S^-,X)$  is obtained by a triangular transformation over $\pi(\Fc_\gamma)$ from   $x_{ik}$, and the system of polynomials $\Rb^+(S^-,X)$ 
is obtained by a triangular transformation over $\pi(\Fc_\gamma)$ from the system of the standard coordinate functions $\{x_{ab}\}$ lying above the anti-diagonal of $X$. This proves 2).~$\Box$\\
\\ 
{\bf Example 4}. For $n=3$ we have
{ \small $$ (Z|X) =  \left(\begin{array}{ccc|ccc}
	z_{11}&z_{12}&z_{13}&	x_{11}&x_{12}&x_{13}\\
	z_{21}&z_{22}&z_{23}&	x_{21}&x_{22}&x_{23}\\
	z_{31}&z_{32}&z_{33}&	x_{31}&x_{32}&x_{33}
	\end{array}
	\right),\qquad  (S^-|X) =  \left(\begin{array}{ccc|ccc}
	0&0&s_{13}&	x_{11}&x_{12}&x_{13}\\
	0&s_{22}&s_{23}&	x_{21}&x_{22}&x_{23}\\
	s_{31}&s_{32}&s_{33}&	x_{31}&x_{32}&x_{33}
	\end{array}
	\right),$$}

{ \small
		$$R^+_{11} (Z,X) = 
	\left|\begin{array}{ccc}
	z_{11}&	z_{12}&x_{11}\\
	z_{21}&	z_{22}&x_{21}\\
	z_{31}&	z_{32}&x_{31}
	\end{array}\right|,\qquad R^+_{11} (S^-,X) = 
	\left|\begin{array}{ccc}
	0&0&x_{11}\\
	0&	s_{22}&x_{21}\\
	s_{31}&	s_{32}&x_{31}
	\end{array}\right| =  D_2(S^-) x_{11}, $$
	$$R^+_{12} (Z,X) = 
 \left|\begin{array}{ccc}
z_{11}&	x_{11}&x_{12}\\
z_{21}&	x_{21}&x_{22}\\
z_{31}&	x_{31}&x_{32}
	\end{array}\right|,\qquad R^+_{12} (S^-,X) = 
	\left|\begin{array}{ccc}
	0&	x_{11}&x_{12}\\
	0&	x_{21}&x_{22}\\
	s_{31}&	x_{31}&x_{32}
	\end{array}\right| =  D_3(S^-) 	\left|\begin{array}{cc}
		x_{11}&x_{12}\\
		x_{21}&x_{22}\\
	\end{array}\right|, $$
	$$R^+_{21} (Z,X) = 
\left|\begin{array}{cc}
z_{21}&	x_{21}\\
z_{31}&	x_{31}
\end{array}\right|,\qquad R^+_{21} (S^-,X) = 
\left|\begin{array}{cc}

0&	x_{21}\\
s_{31}&	x_{31}
\end{array}\right| = - D_3(S^-) x_{21}. $$ }

We define
  $$
\Pb_{\gamma,\tal} (X_\gamma,X_\al)=\left\{\begin{array}{l} \Pb(X_\gamma, X_\al), \mbox{if}~\gamma~\mbox{obeys~ the~ conditions~ of~ Case~~3},\\
\Rb^+(X_\gamma, X_\al), \mbox{if}~\gamma~\mbox{obeys~ the~ conditions~ of~ Case~~4}.
\end{array}\right.
$$
For $h=(X_\al)\in\Hc$, we take
\begin{equation}\label{Pb-tal}
	\Pb_{\tal} (h) = \left\{\begin{array}{l} \varnothing, \quad\mbox{if}\quad \al=\psi(\tal),\\
		 \Pb_{\gamma,\tal} (X_\gamma, X_\al),~ \mbox{where}~ \gamma = \psi(\tal), \quad\mbox{if}\quad \al\ne\psi(\tal). \end{array}\right.
\end{equation}

Let us formulate a statement that follows from Lemmas  \ref{lem-three} and   \ref{lem-four}\\
\Prop\Num\label{prop-tal}. Let $\al$ obeys the conditions of Case 3 or 4, and $\gamma=\psi(\tal)$. Then the system of standard coordinate functions  $\{x_{ab}: a'<b\} $ (they lie above the anti-diagonal of $X_\al$) is algebraically independent on $\Sc_\al$.  
The system of polynomials  $\pi(\Pb_\tal(X_\gamma,X_\al)) $ is obtained by a triangular transformation over $\pi(\Fc)$ from   $\{x_{ab}: a'>b\} $. 

Introduce the notation
\begin{equation}\label{Pb-alpha}
	\Pb_\al=\Pb_\sal\sqcup\Pb_\tal.
\end{equation}

 \Theorem\Num\label{main-theorem}. The field of  $U$-invariants $K(\Hc)^{U}$ is freely generated over $K$ by the system of polynomials
 \begin{equation}\label{system-pd}
 	 \bigcup_{\al\in A}\left(\Pb_{\al}\sqcup \Db_\al\right).
 \end{equation}
  \Proof.  For each $\al\in A$, the system   $\pi( \Db_\al)$ is obtained by a triangular transformation over the field  $K$ from the matrix entries  $\{x_{a,a'}\} $
  that belong to the anti-diagonal of $X_\al$, and $\pi( \Pb_\al)$ is obtained by a triangular transformation over the field  $\pi(\Fc)$  from the other matrix entries of $X_\al$ (see Propositions \ref{prop-sal} and  \ref{prop-tal}).  Therefore, the system of polynomials $\pi(\Pb_\al\sqcup \Db_\al)$ freely generates the field $K(\Sc_\al)$, and 
  $$\bigcup_{\al\in A}\pi(\Pb_\al\sqcup \Db_\al)$$
  freely generates the field $K(\Sc)$. 
 It implies that the embedding $\pi: K(\Hc)^U \hookrightarrow K(\Sc)$ (see Corollary \ref{restriction}) is an isomorphism, and the system (\ref{system-pd}) freely generates $K(\Hc)^{U}$. ~$\Box$
\\ 
\Ex\Num. Let us construct a system of free generators for the quiver  $Q=(V,A)$, defined by the diagram 
$$
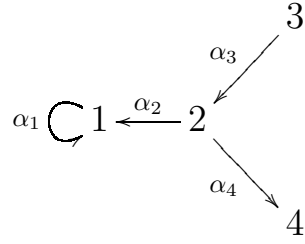
\xymatrix{&&3\ar[ld]_{\al_3}\\
	1\ar@(ul,dl)[]_{\al_1}&2\ar[l]_{\al_2}\ar[rd]_{\al_4}  \\
	&&4}$$

In this case, $V=\{1,2,3,4\}$ and $A=\{\al_1,\al_2,\al_3,\al_4\}$.
The linear space $\Hc$ is a direct sum of four subspaces  $\Hc_i=\Mat(n)$.
The group $GL_Q=\{(g_1,g_2,g_3,g_4)\}$ is a direct product of four copies of  $\GL(n)$.
The action of  $g\in \GL_G$ on $h=(X_1,X_2,X_3,X_3)$ is defined as follows 
$$g.h=(g_1X_1g_1^{-1}, g_1X_2g_2^{-1}, g_2X_3g_3^{-1}, g_4X_4g_2^{-1}).$$
 Let us fix the map $\psi:V\to A$ taking
 $$\psi(1)=\al_1,~~\psi(2)=\psi(4)=\al_4,~~\psi(3) = \al_3.$$
 The choice of $\psi$ defines the section  $\Sc=\oplus_{i=1}^4\Sc_i$ and the systems of polynomials 
 $ \Pb=\sqcup _{i=1}^4\Pb_i$ by the formulas
  $$
\left\{ \begin{array}{l}
	\Sc_1=\Sc^-(n),\\
	\Sc_2=\Mat(n),\\
	\Sc_3=\Sc^+(n),\\
	\Sc_4=\La(n)\\
\end{array}
\right. \quad\mbox{and}\quad \left\{ \begin{array}{l}
	\Pb_1=\Pb(X_1,X_1),\\
	\Pb_2= \Rb^+(X_1,X_2)\sqcup \Rb^-(X_2,X_4),\\
	\Pb_3=\Pb(X_4,X_3),\\
	\Pb_4=\varnothing.\\
\end{array}
\right.
$$
Denote  by $\Db_i=\Db_{\al_i}$ the system of left lower corner minors of $X_i$.

The system of polynomials  $$ \bigsqcup _{i=1}^4~\{\Pb_i\sqcup \Db_i\}$$
freely generates the field  $K(\Hc)^U$.

\end{document}